\documentclass[12pt]{article}
\usepackage{amssymb}
\usepackage{amsmath}
\usepackage{graphicx}
\usepackage{latexsym}
\usepackage{enumerate}
\newtheorem{theorem}{Theorem}[section]
\newtheorem{proposition}[theorem]{Proposition}

\newtheorem{lemma}[theorem]{Lemma}

\newtheorem{remark}[theorem]{Remark}

\newcommand{\bk}{\mathbf{k}}
\newcommand{\bn}{\mathbf{n}}
\newcommand{\ms}{n^\mathbf{r}}
\newcommand{\br}{\mathbf{r}}
\newcommand{\bnp}{\mathbf{n'}}
\newcommand{\perm}{\mathrm{perm}}
\newcommand{\pperm}{\mathrm{pperm}}
\newcommand{\dperm}{\mathrm{dperm}}
\newcommand{\Tr}{\mathrm{Tr}}
\newcommand{\Aut}{\mathrm{Aut}}
\newcommand{\dom}{\mathrm{dom}\,}
\newcommand{\ima}{\mathrm{im}\,}
\newcommand{\indeg}{\mathrm{indeg}\,}
\newcommand{\outdeg}{\mathrm{outdeg}\,}

\begin{document}

\title{Some Generalizations of the MacMahon Master Theorem}
\author{ Michael P. Tuite \and School of Mathematics, Statistics and Applied
Mathematics, \\
National University of Ireland Galway \\
University Road, Galway, Ireland.\\
}
\maketitle

\begin{abstract}
We consider a number of generalizations of the $\beta$-extended MacMahon Master Theorem for a matrix.   
The generalizations are based on replacing permutations on multisets formed from matrix indices by partial permutations or derangements over matrix or submatrix indices.
\end{abstract}

\section{Introduction}
The Master Theorem due to Percy MacMahon first appeared  in 1915 in his classic text \textit{Combinatory Analysis} \cite{MM}.
A generalization known as the $\beta$-extended MacMahon Master Theorem was  discovered in more recent times by Foata and Zeilberger \cite{FZ}.
This present paper is concerned with several further generalizations of the $\beta$-extension informed by recent results in the theory of vertex operator algebras concerning the partition and correlation functions on a genus zero and higher Riemann surface    \cite{MT1,MT2,TZ1,HT,TZ2,TZ3}.

One formulation of the MacMahon Master Theorem (MMT) is the identity of $\det(I-A)^{-1}$, for a given matrix $A$, to an infinite weighted sum over all permanents for matrices indexed by multisets formed from the indices of $A$ \cite{W,KP}. 
The $\beta$-extended MMT relates $\det(I-A)^{-\beta}$ to a similar sum over so-called $\beta$-extended permanents \cite{FZ,KP}.
We consider the following generalizations: 
\begin{enumerate}[(i)]
\item 
The Submatrix MMT. Here the infinite sum runs over multisets formed from the indices of a given submatrix of $A$. 
\item
The Partial Permutation MMT. In this case the $\beta$-extended permanent is replaced by what we refer to as a $(\beta,\theta,\phi)$-extended partial permanent defined in terms of a sum over all partial permutations of the $A$-indices.
\item
The Derangement MMT. We replace the $\beta$-extended permanent by what we refer to as a $\beta$-extended deranged partial permanent defined in terms of a sum over the derangements of the $A$-indices.
\end{enumerate} 

We begin in Section~\ref{Sect_MMT} with a review of the $\beta$-extended MMT \cite{FZ}. 
We provide a graph theoretic proof based on an enumeration of appropriate weights of non-isomorphic permutation graphs labelled by multisets of the indexing set for $A$. In particular, the connected subgraphs are cycles corresponding to permutation cycles.
Section~\ref{Sect_SubMMT} describes our first generalization, the Submatrix MMT 
(Theorem~\ref{Theorem_MMTsub}),   where the set of permutation graphs is modified to account for multisets formed from the indices of an $A$ submatrix. 
In Section~\ref{Sect_PPMMT} we introduce the $(\beta,\theta,\phi)$-extended partial permanent of a matrix, a variation on the $\beta$-extended permanent involving a sum over the partial permutations of the matrix indices. 
The corresponding Partial Permutation MMT (Theorem~\ref{Theorem_MMTpperm}) is proved by a consideration of partial permutation graphs whose connected subgraphs are cycles and open necklaces.  
Section~\ref{Sect_SPPMMT} combines both of the previous generalizations into one general result in Theorem~\ref{MMTsubpperm}.
Finally, in Section~\ref{Sect_DMMT} we introduce another variation, the $\beta$-extended deranged permanent of a matrix, where we sum over the derangements
(fixed point free permutations) of the matrix indices. 
We conclude with a Derangement MMT (Theorem~\ref{Theorem_DerMMT}) and a corresponding Submatrix Derangement MMT (Theorem~\ref{Theorem_dMMTsub})
which are proved by applying the graph theory description of Sections~\ref{Sect_MMT} and \ref{Sect_SubMMT} respectively, but where no 1-cycle graphs occur.

\section{The $\beta$-Extended MacMahon Master Theorem}
\label{Sect_MMT}
Let $A=(A_{ij})$ be an $n\times n$ matrix indexed by $i,j\in \{1,\ldots ,n\}$.
The $\beta$-extended Permanent of $A$ is defined by \cite{FZ}, \cite{KP}%
\begin{equation}
\perm_{\beta }A=\sum\limits_{\pi \in \Sigma_{n}}
\beta ^{C(\pi)}\prod\limits_{i=1}^{n}A_{i\pi (i)},  \label{beta_perm}
\end{equation}
where $C(\pi)$ is the number of cycles in $\pi\in \Sigma_{n}$, the symmetric group.
The  permanent and determinant are the special cases:
\begin{equation}
\perm\, A=\perm_{+1}A,\qquad \det(-A)=\perm_{-1}A.
\label{det_perm}
\end{equation}
Let $\br=(r_{1},\ldots,r_{n})$
denote an $n$-tuple of non-negative integers.
Define
\begin{equation}
\br!=r_{1}!\ldots r_{n}!,
\label{rfact}
 \end{equation} 
and let
\begin{equation}
\ms=\{1^{r_{1}}2^{r_{2}}\ldots n^{r_{n}}\} 
=\{1_{1},\ldots,1_{r_{1}},\ldots, n_{1},\ldots,n_{r_{n}}\}, \label{mset}
\end{equation}
denote the multiset of size $N=\sum\limits_{i=1}^{n}r_{i} $ formed from the original index set 
$\{1,\ldots,n\} $ where the index $i$ is repeated $r_{i}$ times. 
We sometimes notate a repeated index by $i_{a}$ for label $a=1,\ldots, r_{i}$. 
For an $n\times n$ matrix $A$, we let $A(\ms,\ms)$ denote the $N\times N$ matrix indexed by the elements of $\ms$ and define $A(\ms,\ms)=1$ for $\br=(0,0,\ldots,0)$.

We now describe a generalization,  due to Foata and Zeilberger \cite{FZ}, of the  MacMahon Master Theorem
(MMT) of classical combinatorics \cite{MM}. 
We give a detailed proof based on a graph theory method which is extensively employed throughout this paper. This proof is very similar to that of Theorem~5 of \cite{MT1} where the MMT was essentially rediscovered.  
\begin{theorem}[The $\beta$-Extended MMT]
\label{Theorem_MMT} 
\begin{equation}
\sum_{r_{i}\ge 0}
\frac{1}{\br!} \,\perm_{\beta} A(\ms,\ms)
=\frac{1}{\det (I-A)^{\beta }}.  \label{MMT}
\end{equation}
\end{theorem}

\begin{remark}
For $\beta=1$, Theorem~\ref{Theorem_MMT} reduces to the MMT. For $\beta=-1$ we use \eqref{det_perm} to find that only proper subsets of $\{1,\ldots ,n\}$ contribute resulting in the
determinant identity for $B=-A$ e.g. \cite{TZ1}
\begin{equation*}
\sum\limits_{r_{i}\in\{ 0,1\}}\det B(\ms,\ms)=\det (I+B).
\end{equation*}
\end{remark}

\textbf{Proof of Theorem~\ref{Theorem_MMT}.} 
Let $ \Sigma (\ms)$ denote the symmetric group of the multiset $\ms$. 
For $\pi\in\Sigma(\ms)$ we define a permutation graph $\gamma_{\pi}$
with $N$ vertices labelled by $i\in\{1,\ldots,n\}$, 
and with directed edges 
\begin{equation*}
e_{ij}=i\,\bullet \longrightarrow \bullet \,j\,,
\end{equation*}%
provided $j=\pi (i)$.  
The connected subgraphs of $\gamma_{\pi}\in \Gamma$ are cycles arising from the cycles of $\pi$.
For example, for $n=4$ with $\br=(3,2,0,1)$ and permutation 
$\pi=(1_{1}2_{1}1_{2}2_{2})(1_{3}4_{1})$ the corresponding graph has two cycles as shown in Fig.~1
\begin{center}
\begin{picture}(280,70)

\put(100,50){\line(1,0){40}}
\put(100,50){\vector(1,0){23}}
\put(90,50){\makebox(1,0){$1$}}
\put(100,50){\circle*{4}}

\put(140,50){\line(0,-1){40}}
\put(140,50){\vector(0,-1){23}}
\put(150,50){\makebox(1,0){$2$}}
\put(140,50){\circle*{4}}

\put(140,10){\line(-1,0){40}}
\put(140,10){\vector(-1,0){23}}
\put(150,10){\makebox(1,0){$1$}}
\put(140,10){\circle*{4}}

\put(100,10){\line(0,1){40}}
\put(100,10){\vector(0,1){23}}
\put(90,10){\makebox(1,0){$2$}}
\put(100,10){\circle*{4}}

\put(180,30){\qbezier(0,0)(20,20)(40,0)}
\put(200,40){\vector(1,0){2}}
\put(170,30){\makebox(1,0){$1$}}
\put(180,30){\circle*{4}}

\put(180,30){\qbezier(0,0)(20,-20)(40,0)}
\put(200,20){\vector(-1,0){2}}
\put(230,30){\makebox(1,0){$4$}}
\put(220,30){\circle*{4}}
\end{picture}
\newline
{\small Fig.~1 $\gamma_{\pi}$ for $\pi=(1_{1}2_{1}1_{2}2_{2})(1_{3}4_{1})$. }
\end{center}

Define a weight for each edge of $\gamma_{\pi }$ by
\begin{eqnarray*} 
w(e_{ij})=A_{ij},
\end{eqnarray*}
and a weight for $\gamma_{\pi }$ by 
\begin{eqnarray}
w(\gamma_{\pi })&=&
\beta^{C(\pi )}
\prod_{e_{ij}\in \gamma_{\pi }}w(e_{ij}).\label{wtgam}
\end{eqnarray}%
where $C(\pi )$ is the number of cycles in $\pi $. 
Note that the weight is multiplicative with respect to the cycle decomposition of $\pi$.
\eqref{wtgam} also implies
\begin{equation}
\perm_{\beta } A(\ms,\ms)
=\sum_{\pi \in \Sigma (\ms)}w(\gamma_{\pi }).
\label{permA}
\end{equation}%

Let $\Lambda(\br)=\Sigma_{r_{1}}\times \ldots \times\Sigma_{r_{n}}\subseteq\Sigma (\ms)$ denote the label group of order $\vert\Lambda(\br)\vert=\br!$ which permutes the identical elements of $\ms$.
$\Lambda(\br)$ generates isomorphic graphs with $\gamma_{\pi}\sim\gamma_{\lambda \pi \lambda^{-1}}$ for $\lambda\in\Lambda(\br)$
and the automorphism group of $\gamma_{\pi}$ is the $\pi $ stabilizer  $\Aut(\gamma_{\pi})=\{\lambda \in \Lambda(\br)\vert \lambda\pi=\pi\lambda \}\subseteq\Lambda(\br)$.  
Using the Orbit-Stabilizer theorem it follows that the number of isomorphic graphs generated by the action of $\Lambda(\br)$ on $\gamma_{\pi}$ is given by
\begin{equation}
\vert\Lambda(\br)\gamma_{\pi}\vert = \frac{\vert\Lambda(\br)\vert}{\vert\Aut(\gamma_{\pi})\vert}. 
\label{Lambdag}
\end{equation}  
(e.g. in Fig.~1,  $\Lambda(\br)=\Sigma_{2}\times \Sigma_{3}$ and $\Aut(\gamma_{\pi})=\Sigma_{2}$ so that there are 6 permutations in $\Sigma(\ms)$ with graph $\gamma_{\pi}$). 
Combining \eqref{permA} and  \eqref{Lambdag} we find that
\begin{eqnarray}
\sum_{\br}\frac{1}{\br!}\,\perm_{\beta } A(\ms,\ms) 
&=& 
\sum_{\gamma\in \Gamma} \frac{w(\gamma)}{\vert \Aut(\gamma) \vert},
\label{gammsum}
\end{eqnarray}
where $\Gamma$ denotes the set of non-isomorphic graphs.

Consider the decomposition of a graph $\gamma$ into cycle graphs
\begin{equation*}
\gamma =\gamma_{\sigma_{1}}^{m_{1}}\ldots \gamma_{\sigma_{K}}^{m_{K}}, 
\end{equation*}%
where $\{\gamma_{\sigma_{i}}\}$ are non-isomorphic and $\gamma_{\sigma_{i}}$ occurs $m_{i}$ times. 
The automorphism group is 
\begin{equation*}
\Aut(\gamma)=\prod_{i=1}^{M}\Aut(\gamma_{\sigma_{i}}^{m_{i}}),
\end{equation*}
where $\Aut(\gamma_{\sigma}^{m})=\Sigma_{m}\rtimes \Aut(\gamma_{\sigma})^{m}$
of order $m!\,|\Aut(\gamma_{\sigma})|^{m}$. 
Furthermore, since the weight is multiplicative,
$
w (\gamma )=\prod_{i=1}^{M}w (\gamma_{\sigma_{i}})^{m_{i}}$.
Thus we find
\begin{eqnarray}
\sum_{\gamma \in \Gamma }\frac{w (\gamma )}{|\Aut(\gamma )|}
&=&\prod_{\gamma_{\sigma }\in \Gamma_{\sigma }}\sum_{m\geq 0}\frac{1}{m!}
\left(
\frac{w (\gamma_{\sigma })}{|\Aut(\gamma_{\sigma })|}
\right)^{m}
\notag \\
&=&\exp \left( \sum_{\gamma_{\sigma }\in \Gamma_{\sigma }}\frac{w(\gamma_{\sigma })}{|\Aut(\gamma_{\sigma })|}\right),
\label{expwt}
\end{eqnarray}
where $\Gamma_{\sigma }$ denotes the set of non-isomorphic cycle graphs.
For a cycle $\sigma $ of order $|\sigma|=t$ we have 
$\Aut(\gamma_{\sigma })=\langle \sigma ^{s}\rangle $ for some $s|t$ with $|\Aut(\gamma_{\sigma })|=\frac{t}{s}$.
Using the trace identity 
\begin{equation*}
\sum_{\gamma_{\sigma },|\sigma |=t}s\ w (\gamma_{\sigma })
=\beta \Tr(A^{t}),
\end{equation*}%
we find 
\begin{eqnarray*}
\sum_{\gamma_{\sigma }\in \Gamma_{\sigma }}\frac{w (\gamma
_{\sigma })}{|\Aut(\gamma_{\sigma })|} &=&\beta \sum_{t\geq 1}\frac{%
1}{t}\Tr(A^{t}) \\
&=&-\beta \mathrm{Tr}\log (I-A) \\
&=&-\beta \log \det (I-A).
\end{eqnarray*}%
Thus 
\begin{equation*}
\sum_{\br}\frac{1}{\br!}\,\perm_{\beta } A(\ms,\ms)=\det (I-A)^{-\beta }.\qquad\square 
\end{equation*}
\bigskip

Let $w_{1}(\gamma)$ denote the weight for $\gamma$ with $\beta =1$ in \eqref{wtgam}. 
Define a cycle to be primitive (or rotationless) 
if $|\Aut(\gamma_{\sigma })|=1$. For a general cycle $\sigma $ with 
$|\Aut(\gamma_{\sigma })|=k$ we have
$\gamma_{\sigma }=\gamma_{\rho }^{k}$ 
for a primitive cycle $\rho $. Let $\Gamma_{\rho }$ denote the set of
all primitive cycles. Then 
\begin{eqnarray*}
\sum_{\gamma_{\sigma }\in \Gamma_{\sigma }}
\frac{w_{1}(\gamma_{\sigma })}{|\Aut(\gamma_{\sigma })|} 
&=&\sum_{\gamma_{\rho }\in 
\Gamma_{\rho}}
\sum_{k\geq 1}\frac{1}{k}w_{1}(\gamma_{\rho })^{k} \\
&=&-\sum_{\gamma_{\rho }\in \Gamma_{\rho }}\log \det (1-w_{1}(\gamma
_{\rho })).
\end{eqnarray*}%
Combining this with \eqref{expwt} implies \cite{MT1}

\begin{proposition}
\begin{equation*}
\det (I-A)=\prod_{\gamma_{\rho }\in \Gamma_{\rho }}(1-w_{1}(\gamma_{\rho
})).
\end{equation*}
\end{proposition}

\bigskip
 
\section{The Submatrix MMT}
\label{Sect_SubMMT}
Our first generalization of Theorem \ref{Theorem_MMT} concerns submatrices. 
Consider an $(n^{\prime }+n)\times (n^{\prime }+n)$
matrix with block structure 
\begin{equation}
\left[
\begin{array}{cc}
B & U \\ 
V & A
\end{array}
\right],
\label{Block}
\end{equation}%
where $A=(A_{ij})$ is an $n\times n$ matrix indexed by $i,j$, $B=(B_{i^{\prime }j^{\prime }})$ is an $n^{\prime }\times n^{\prime }$
matrix indexed by $i^{\prime },j^{\prime }$, $U=(U_{i^{\prime }j})$ is an $n^{\prime }\times n$ matrix
and $V=(V_{ij^{\prime }})$ is an $n\times n^{\prime }$ matrix. For a
multiset $\ms$ of size $N$ define 
 the $(n^{\prime }+N)\times (n^{\prime }+N)$ matrix%
\begin{equation}
\left[
\begin{array}{cc}
B & U(\ms) \\ 
V(\ms) & A(\ms,\ms)
\end{array}
\right],
\label{BAk}
\end{equation}%
where, as before, $A(\ms,\ms)$ denotes the $N\times N$ matrix indexed by 
$\ms$, $U(\ms)$ is an $n^{\prime }\times N$ matrix and $V(\ms)$ is an 
$N\times n^{\prime }$ matrix. We then find

\begin{theorem}
\label{Theorem_MMTsub} 
\begin{eqnarray}
\sum_{\br}\frac{1}
{\br!}\,\perm_{\beta }
\left[\begin{array}{cc}
B & U(\ms) \\ 
V(\ms) & A(\ms,\ms)
\end{array}\right] 
&=&\frac{\perm_{\beta }\widetilde{B}}{\det
(I-A)^{\beta }},
\label{MMTsub}
\end{eqnarray}%
for $n'\times n'$ matrix  
\begin{equation*}
\widetilde{B}=B+U(I-A)^{-1}V,
\end{equation*}
where $(I-A)^{-1}=\sum_{k\ge 0}A^{k}$.
\end{theorem}
This result is related to Theorem~10 of \cite{MT1} when $\beta=1$ and Theorem~2 of \cite{TZ1} for $\beta=-1$.
\bigskip

\noindent \textbf{Proof.} 
Let $\bn=\{1,\ldots ,n\}$ and 
$\bnp=\{1^{\prime },\ldots ,n^{\prime }\}$ and let  
$\bnp\cup \ms$ denote the multiset indexing the block matrix \eqref{BAk}.
Define a permutation graph $\gamma_{\pi }$ with weight $w (\gamma_{\pi })$ for each $\pi \in \Sigma (\bnp\cup \ms)$ as follows.
Each vertex is labelled by an element of $\bn$ or  $\bnp$ which we refer to as $\bn$-vertex or $\bnp$-vertex respectively. 
For $	l= \pi(k)$ with $k,l\in \bnp\cup \ms$ we define an edge $e_{kl}=k\, \bullet \longrightarrow {\bullet }\,l$ with weight
\begin{equation*}
 w(e_{kl})=\left[
\begin{array}{cc}
B & U(\ms) \\ 
V(\ms) & A(\ms,\ms)
\end{array}
\right]_{kl}. 
\end{equation*}%
Define a weight for $\gamma_{\pi}$ by
\begin{eqnarray*} 
w (\gamma_{\pi }) &=&\beta ^{C(\pi )}
\prod_{e_{kl}\in \gamma_{\pi}} w (e_{kl}),
\end{eqnarray*}%
where $C(\pi )$ is the number of cycles in $\pi $. 
As before, we find 
\begin{equation*}
\sum_{\br}\frac{1}{\br!}\,
\perm_{\beta } \left[\begin{array}{cc}
B & U(\bk) \\ 
V(\bk) & A(\bk,\bk)
\end{array}\right]=\sum_{\gamma \in \widehat{\Gamma}}\frac{w(\gamma )}{|\Aut(\gamma )|},
\end{equation*}%
where $\widehat{\Gamma}$ denotes the set of non-isomorphic graphs.
Each $\gamma\in \widehat{\Gamma}$ has a decomposition into cycles $\gamma_{\sigma_{a}}$ 
which contain $\bn$-vertices only and cycles $\gamma_{\sigma^{\prime}_{b}}$ which
contain at least one $\mathbf{n'}$-vertex:  
\begin{equation*}
\gamma =\gamma_{\sigma_{1}}^{m_{1}}\ldots
\gamma_{\sigma_{K}}^{m_{K}}
\gamma_{\sigma^{\prime} _{1}}\ldots \gamma_{\sigma^{\prime} _{L}},
\end{equation*}%
with weight 
\begin{equation*}
w (\gamma)=
\prod_{a}w (\gamma_{\sigma_{a}})^{m_{a}}
\prod_{b}w (\gamma_{\sigma^{\prime}_{b}}).
\end{equation*}%
The set of non-isomorphic $\gamma_{\sigma_{a}}$ cycle graphs labelled by $\bn$ is equivalent to $\Gamma_{\sigma }$ introduced in the proof of Theorem~\ref{Theorem_MMT}.
Since each $\mathbf{n'}$-vertex occurs exactly once in $\gamma$, each $\gamma_{\sigma^{\prime}_{b}}$ cycle 
occurs at most once and has trivial automorphism group. Hence 
\begin{equation*}
|\Aut(\gamma)|=
\prod_{a}|\Aut(\gamma_{\sigma_{a}})|^{m_{a}}m_{a}!,
\end{equation*}%
as before.
Thus the sum over weights of all graphs decomposes into the product%
\begin{eqnarray*}
\sum_{\gamma \in \widehat{\Gamma}}\frac{w (\gamma )}{|\Aut(\gamma)|}
&=&
\sum_{\gamma_{\sigma^{\prime}} }w (\gamma_{\sigma^{\prime}} )
\prod_{\gamma_{\sigma }\in \Gamma_{\sigma }}\sum_{m\geq 0}\frac{%
w (\gamma_{\sigma })^{m}}{|\Aut(\gamma_{\sigma })|^{m}m!}\\
&=&\frac{\sum_{\gamma_{\sigma^{\prime}}}w (\gamma_{\sigma^{\prime}})}{\det(I-A)^{\beta }},
\end{eqnarray*}%
using Theorem~\ref{Theorem_MMT} and where $\gamma_{\sigma^{\prime}}$ ranges over non-isomorphic cycles in $\widehat{\Gamma}$ containing at least one $\mathbf{n'}$-vertex.

It remains to compute $\sum_{\gamma_{\sigma^{\prime}} }w (\gamma_{\sigma^{\prime}} )$. 
Let $\sigma'\in\Sigma_{n'}$ denote the permutation cycle corresponding to the cyclic sequence of $\mathbf{n'}$-vertices in a $\gamma_{\sigma^{\prime}}$-cycle (for arbitrary intermediate $\bn$-vertices). 
The total edge weight coming from all subgraphs, illustrated in Fig.~2, joining two $\mathbf{n'}$-vertices, 
$i^{\prime }$ and $j^{\prime}$, 
summed over all intermediate $\bn$-vertices is
\begin{eqnarray*}
&&B_{i^{\prime }j^{\prime }}
+(UV)_{i^{\prime}j^{\prime }}
+(UAV)_{i^{\prime}j^{\prime }}+(UA^2V)_{i^{\prime}j^{\prime }}+\ldots  \\
&=&(B+U(I-A)^{-1}V)_{i^{\prime }j^{\prime }}=\widetilde{B}_{i^{\prime }j^{\prime }}.
\end{eqnarray*}%
\begin{center}
\begin{picture}(340,40)

\put(10,10){\line(1,0){40}}
\put(10,10){\vector(1,0){23}}
\put(10,20){\makebox(1,0){$i^{\prime}$}}
\put(10,10){\circle*{4}}

\put(50,20){\makebox(1,0){$j'$}}
\put(50,10){\circle*{4}}
\put(60,10){\makebox(1,0){,}}


\put(80,10){\line(1,0){40}}
\put(80,10){\vector(1,0){23}}
\put(80,20){\makebox(1,0){$i^{\prime}$}}
\put(80,10){\circle*{4}}

\put(120,10){\line(1,0){40}}
\put(120,10){\vector(1,0){23}}
\put(120,20){\makebox(1,0){$k$}}
\put(120,10){\circle*{4}}

\put(160,20){\makebox(1,0){$j'$}}
\put(160,10){\circle*{4}}

\put(170,10){\makebox(1,0){,}}


\put(190,10){\line(1,0){40}}
\put(190,10){\vector(1,0){23}}
\put(190,20){\makebox(1,0){$i^{\prime}$}}
\put(190,10){\circle*{4}}

\put(230,10){\line(1,0){40}}
\put(230,10){\vector(1,0){23}}
\put(230,20){\makebox(1,0){$k$}}
\put(230,10){\circle*{4}}

\put(270,10){\line(1,0){40}}
\put(270,10){\vector(1,0){23}}
\put(270,20){\makebox(1,0){$l$}}
\put(270,10){\circle*{4}}

\put(310,20){\makebox(1,0){$j'$}}
\put(310,10){\circle*{4}}

\put(330,10){\makebox(1,0){$,\ldots$}}
\end{picture}

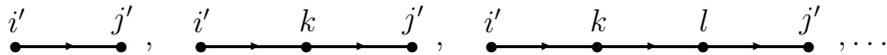
{\small Fig.~2 }
\end{center}
Thus the total weight of all $\gamma_{\sigma^{\prime}}$ cycles for a given $\mathbf{n'}$-vertex cycle $\sigma'=(i'_{1}\ldots i'_{p})$ is $\beta \prod_{l} \widetilde{B}_{i'_{l}\sigma'(i'_{l})}$. Altogether, it follows that 
\begin{eqnarray*}
\sum_{\gamma_{\sigma^{\prime}}}w (\gamma_{\sigma^{\prime}} ) &=&\sum\limits_{\pi ^{\prime }\in \Sigma
_{n^{\prime }}}\beta ^{C(\pi ^{\prime })}
\prod\limits_{i^{\prime}}
\widetilde{B}_{i^{\prime }\pi ^{\prime }(i^{\prime })} \\
&=&\perm_{\beta }\widetilde{B}.\quad \square
\end{eqnarray*}
\medskip

\begin{lemma}
For $\beta=-1$ Theorem~\ref{Theorem_MMTsub} implies
\begin{eqnarray*}
\sum_{r_{i}\in\{ 0,1\}}\det
\left[\begin{array}{cc}
B & U(\ms) \\ 
V(\ms) & A(\ms,\ms)
\end{array}\right] 
&=&\det
\left[\begin{array}{cc}
B & -U \\ 
-V & I+A
\end{array}\right].
\end{eqnarray*}
\end{lemma}

\noindent \textbf{Proof.} For $\beta=-1$ the right hand side of \eqref{MMTsub} gives
\begin{eqnarray*}
\perm_{-1}\widetilde{B}\,\det(I-A)&=&(-1)^{n^{\prime}}\det(B+U(I-A)^{-1}V)\det(I-A)\\
&=&
\det\left[ 
\begin{array}{cc}
-B & U   
\\
V & I - A
\end{array}
\right], 
\end{eqnarray*}
by means of the matrix identity 
\begin{eqnarray*}
&
\left[ 
\begin{array}{cc}
-B & U   
\\
V & I - A
\end{array}
\right] 
 = 
 \left[ 
\begin{array}{cc}
-I^{\prime} & U(I - A)^{-1}  
\\
0 & I 
\end{array}
\right]
\left[ 
\begin{array}{cc}
B+U(I - A)^{-1}V  & 0  
\\
V & I 
\end{array}
\right]
\left[ 
\begin{array}{cc}
I^{\prime} & 0   
\\
0 & I - A
\end{array}
\right],&
\end{eqnarray*}
where $I$ and $I^{\prime}$ are respectively $n\times n$ and $n^{\prime}\times n^{\prime}$ identity matrices. 
The result follows on replacing $A,B,U,V$ by $-A,-B,-U,-V$.\quad $\square$ 

\bigskip 
\section{The Partial Permutation MMT}
\label{Sect_PPMMT}
The next generalization of Theorem~\ref{Theorem_MMT} is concerned with replacing  
permutations by partial permutations with a suitable
 generalization of the notions of permanent and $\beta$-extended permanent.
Let $\Psi$ denote the set of partial permutations of the set  $\{1,\ldots ,n\}$ 
i.e. injective partial mappings from $\{1,\ldots ,n\}$ to itself. 
For $\psi\in\Psi$ we let $\dom \psi$ and $\ima \psi$ denote the domain and image respectively and
let $\pi_{\psi}$ denote the (possibly empty) permutation of $\dom \psi \cap \ima \psi$ determined by $\psi$.

We introduce the Partial Permanent of an $n\times n$ matrix $A=(A_{ij})$ indexed by $i,j\in \{1,\ldots ,n\}$ as follows 
\begin{equation}
\pperm A= \sum\limits_{\psi\in\Psi }\prod\limits_{i \in\, \dom \psi} A_{i\psi(i)},
\label{pperm}
\end{equation}
with unit contribution for the empty map. 
Let $\theta=(\theta_{i}),\phi=(\phi_{i})$ be $n$-vectors and define the $(\beta,\theta,\phi)$-extended Partial Permanent by 
\begin{equation}
\pperm_{\beta\theta\phi} A= \sum\limits_{\psi\in\Psi }
\beta^{C(\pi_{\psi})}
\prod\limits_{i \in \dom \psi} A_{i\psi(i)}
\prod\limits_{j \not\in \, \ima \psi} \theta_{j}
\prod\limits_{k \not\in \, \dom \psi} \phi_{k}
,
\label{ppermBTP}
\end{equation}
where $C(\pi_{\psi})$ is the number of cycles in $\pi_{\psi}$ e.g. 
\begin{eqnarray*}
\pperm_{\beta\theta\phi}\left[ 
\begin{array}{cc}
A_{11} & A_{12}\\
A_{21} & A_{22}	
\end{array}
\right]
&=&\theta_{1}\phi_{1}\theta_{2}\phi_{2}+\beta(A_{11}\theta_{2}\phi_{2}+A_{22}\theta_{1}\phi_{1})\\
&&+A_{12}\theta_{1}\phi_{2}+A_{21}\theta_{2}\phi_{1}+\beta^2 A_{11}A_{22}+\beta A_{12}A_{21}.
\end{eqnarray*} 
A recent application of an extended partial permanent appears in
\cite{HT}.
\medskip

Let $A(\ms,\ms)$ denote the $N\times N$ matrix indexed by a multiset
$\ms$ as before. 
We also let  $\pperm_{\beta\theta\phi}A(\ms,\ms)$ denote the corresponding partial permanent with $N$-vectors  $(\theta_{1_1},\ldots,\theta_{n_{r_{n}}})$ and $(\phi_{1_1},\ldots,\phi_{n_{r_{n}}})$. 
We then find
\begin{theorem}
\label{Theorem_MMTpperm} 
\begin{eqnarray}
\sum_{\br}\frac{1}{\br!}\, \pperm_{\beta\theta\phi }
A(\ms,\ms)
&=&\frac{
e^{\theta(I-A)^{-1}\phi^{T}}
}
{\det (I-A)^{\beta }},
\label{ppermsum}
\end{eqnarray}%
where $\phi^{T}$ denotes the transpose of the row vector $\phi$.
\end{theorem}
This result is related to Theorem~11 of \cite{MT1} for $\beta=1$.
\bigskip

\noindent \textbf{Proof.} 
Let $\Psi(\ms)$ denote the partial permutations of $\ms$. 
Define a partial permutation graph $\gamma_{\psi }$ labelled by $\{1,\ldots ,n\}$ for each $\psi \in \Psi(\bk)$ with 
edges
\begin{equation*}
e_{ij}=i\, \bullet \longrightarrow \bullet\, j\, ,
\end{equation*}%
for $j=\psi (i)$ with $i\in\dom\psi$ and $j\in \ima\psi$.  
Let $v_{i}$ denote the vertex of $\gamma_{\psi }$ with label $i$. 
If $i\notin \dom\psi$ then either $\deg v_{i}= 0$ or $\deg v_{i}= \indeg v_{k}=1$
whereas if $i\notin \ima\psi$ then either $\deg v_{i}= 0$ or
$\deg v_{i}= \outdeg v_{i}=1$. 
In all other cases $\deg v_{i}= 2$ with $ \indeg v_{i}=\outdeg v_{i}=1$. 
The connected subgraphs in this case consist of cycles and  open necklaces i.e. graphs with two end points of degree one. 
We regard a graph consisting of a single degree zero vertex as a degenerate necklace.
For example, for $n=4$, $\br=(3,2,0,1)$ and partial permutation 
$\psi=
\left(
\begin{array}{cccccc}
1_{1} & 1_{2} & 1_{3} & 2_{1} & 2_{2} & 4_{1}\\
2_{1} &       & 4_{1} & 1_{2} &       & 1_{3}
\end{array}
\right)$ then $\gamma_{\psi}$ is shown in Fig.~3. In this case 
$\dom \psi=\{1_{1}, 1_{3}, 2_{1},4_{1}\}$ and $\ima \psi=\{ 1_{2},1_{3},2_{1} ,4_{1}\}$ 
and $\pi_{\psi}=(1_{3}4_{1})$.
\begin{center}
\begin{picture}(280,60)

\put(100,50){\line(1,0){40}}
\put(100,50){\vector(1,0){23}}
\put(90,50){\makebox(1,0){1}}
\put(100,50){\circle*{4}}

\put(140,50){\line(0,-1){40}}
\put(140,50){\vector(0,-1){23}}
\put(150,50){\makebox(1,0){2}}
\put(140,50){\circle*{4}}

\put(150,10){\makebox(1,0){1}}
\put(140,10){\circle*{4}}

\put(90,10){\makebox(1,0){2}}
\put(100,10){\circle*{4}}

\put(180,30){\qbezier(0,0)(20,20)(40,0)}
\put(200,40){\vector(1,0){2}}
\put(170,30){\makebox(1,0){1}}
\put(180,30){\circle*{4}}

\put(180,30){\qbezier(0,0)(20,-20)(40,0)}
\put(200,20){\vector(-1,0){2}}
\put(230,30){\makebox(1,0){4}}
\put(220,30){\circle*{4}}
\end{picture}

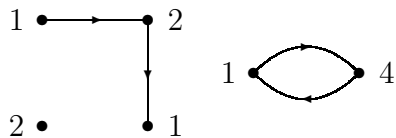
{\small Fig.~3 }
\end{center}

Define an edge weight as before by $w (e_{ij}) = A_{ij}$ and introduce a vertex weight   
\begin{equation*}
w(v_{k})=
\left\{
\begin{array}{cl}
1, & 	\deg v_{k}= 2,\\
\theta_{k}, & 	\deg v_{k}= \outdeg v_{k}=1,\\
\phi_{k}, & 	\deg v_{k}= \indeg v_{k}=1,\\
\theta_{k}\phi_{k}, & 	\deg v_{k}= 0.
\end{array}
\right. 
\end{equation*}
The weight of a graph $\gamma_{\psi}$ is defined by
\begin{eqnarray*}
w (\gamma_{\psi}) &=&\beta ^{C(\pi_{\psi} )}
\prod_{e_{ij}} w (e_{ij})
\prod_{v_{k}} w(v_{k}),
\end{eqnarray*}%
where $C(\pi_{\psi} )$ is the number of cycles in $\pi_{\psi} $. 
The weight is multiplicative with respect to the cycle and necklace decomposition.
We find again that
\begin{equation*}
\sum_{\br}\frac{1}{\br!}\pperm_{\beta\theta\phi } \, A(\ms,\ms)
=\sum_{\gamma \in \widetilde{\Gamma}}\frac{w(\gamma )}{|\Aut(\gamma )|},
\end{equation*}%
where $\widetilde{\Gamma}$ denotes the set of non-isomorphic graphs.
Each $\gamma\in \widetilde{\Gamma}$ has a decomposition into connected  cycle graphs $\gamma_{\sigma_{a}}$ and open necklaces $\nu_{b}$:  
\begin{equation*}
\gamma =\nu_{1}^{l_{1}}\ldots \nu_{L}^{l_{1}}
\gamma_{\sigma_{1}}^{m_{1}}\ldots
\gamma_{\sigma_{K}}^{m_{K}},
\end{equation*}%
with weight 
\begin{equation*}
w (\gamma)=
\prod_{b}w (\nu_{b})^{l_{b}}\cdot
\prod_{a}w (\gamma_{\sigma_{a}})^{m_{a}}.
\end{equation*}%
Each necklace has trivial automorphism group but can have multiple occurrences. 
Hence we find that
\begin{equation*}
|\Aut(\gamma)|=\prod_{b}l_{b}!
\cdot
\prod_{a}|\Aut(\gamma_{\sigma_{a}})|^{m_{a}}m_{a}!.
\end{equation*}%
Thus the sum over weights of all graphs decomposes into the product%
\begin{eqnarray*}
\sum_{\gamma \in\, \widetilde{\Gamma}}\frac{w (\gamma )}{|\Aut(\gamma)|}
&=&
\prod_{\nu\in\,\Gamma_{\nu}}
\sum_{l\geq 0}\frac{w (\nu )^{l}}{l!}
\cdot
\prod_{\gamma_{\sigma }\in\, \Gamma_{\sigma }}
\sum_{m\geq 0}\frac{w (\gamma_{\sigma })^{m}}{|\Aut(\gamma_{\sigma })|^{m}m!}
\\
&=&\exp\left(\sum_{\nu\in\,\Gamma_{\nu} }w (\nu )\right)\frac{1}
{\det(I-A)^{\beta }},
\end{eqnarray*}%
where $\Gamma_{\nu}$ denotes the set of non-isomorphic open necklaces
and using Theorem~\ref{Theorem_MMT} again.
Finally,  the sum over the weights of connected necklaces, such as depicted in Fig.~4, is
\begin{eqnarray*}
\sum_{\nu \in\, \Gamma_{\nu}}w (\nu )&=&
\theta\phi^T
+\theta A \phi^T
+\theta A^2 \phi^T+\ldots  \\
&=&\theta(I-A)^{-1}\phi^T.\quad \square
\end{eqnarray*}%
\begin{center}
\begin{picture}(340,40)


\put(50,20){\makebox(1,0){$i$}}
\put(50,10){\circle*{4}}
\put(60,10){\makebox(1,0){,}}


\put(80,10){\line(1,0){40}}
\put(80,10){\vector(1,0){23}}
\put(80,20){\makebox(1,0){$i$}}
\put(80,10){\circle*{4}}

\put(120,20){\makebox(1,0){$j$}}
\put(120,10){\circle*{4}}

\put(130,10){\makebox(1,0){,}}


\put(150,10){\line(1,0){40}}
\put(150,10){\vector(1,0){23}}
\put(150,20){\makebox(1,0){$i$}}
\put(150,10){\circle*{4}}

\put(190,10){\line(1,0){40}}
\put(190,10){\vector(1,0){23}}
\put(190,20){\makebox(1,0){$j$}}
\put(190,10){\circle*{4}}

\put(230,20){\makebox(1,0){$k$}}
\put(230,10){\circle*{4}}

\put(250,10){\makebox(1,0){$,\ldots$}}

\end{picture}

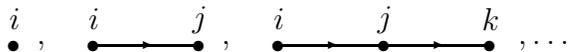
{\small Fig.~4 }
\end{center}

\bigskip
\textbf{Example.} Consider $n=1$ with $A=z$ and $\theta_{1}=\phi_{1}=\sqrt{\alpha z}$. Then we find 
\begin{equation*}
\pperm_{\beta\theta\phi }A(1^{r},1^{r})=p_{r}(\alpha,\beta)z^{r},
\end{equation*}
where $p_{r}(\alpha,\beta)=\sum_{s,t}p_{rst}\alpha^{s}\beta^{t}$ is the generating polynomial for $p_{rst}$ the number of graphs with $r$ identically labelled vertices, $s$ open necklaces and $t$ cycles.  Theorem~\ref{Theorem_MMTpperm} provides the exponential generating function for $p_{r}(\alpha,\beta)$ \cite{HT}
\begin{equation*}
\sum_{r\ge 0} \frac{p_{r}(\alpha,\beta)}{r!}z^{r}=
\frac{\exp\left( \frac{\alpha z}{1-z}\right)}{(1-z)^\beta}.
\end{equation*}

\bigskip 
\section{The Submatrix Partial Permutation MMT}
\label{Sect_SPPMMT}
We can combine the two generalizations above into one theorem concerning partial permutations of submatrices of the $(n'+n)\times (n'+n)$ block matrix  \eqref{Block}. 
Let $\theta'=(\theta'_{i'})$ and $\phi'=(\phi'_{i'})$ 
be $n'$-vectors and $\theta=(\theta_{i})$ and
 $\phi=(\phi_{i})$ be $n$-vectors. 
 For a multiset $\ms$ of size $N$ and block matrix \eqref{BAk} labelled by 
$\mathbf{n'}=\{1',\ldots,n'\}$ and $\ms$, we let 
$\pperm_{\beta \theta\phi}
\left[\begin{array}{cc}
B & U(\ms) \\ 
V(\ms) & A(\ms,\ms)
\end{array}\right]$ 
denote the $(\beta,\theta,\phi)$-extended partial permanent with $(n'+N)$-vectors 
$(\theta'_{1'},\ldots,\theta'_{n'},\theta_{1_{1}},\ldots,\theta_{n_{r_{n}}})$ and 
$(\phi'_{1'},\ldots,\phi'_{n'},\phi_{1_{1}},\ldots,\phi_{n_{r_{n}}})$ respectively. 
We then find
\begin{theorem}
\label{MMTsubpperm} 
\begin{eqnarray}
\sum_{\br}\frac{1}{\br!}\, \pperm_{\beta \theta\phi}
\left[\begin{array}{cc}
B & U(\ms) \\ 
V(\ms) & A(\ms,\ms)
\end{array}\right] 
&=&
\frac{
e^{\theta(I-A)^{-1}\phi^{T}}\cdot
\pperm_{\beta\,\widetilde{\theta}\,\widetilde{\phi}}\widetilde{B}
}
{\det (I-A)^{\beta }},\quad
\label{subppermsum}
\end{eqnarray}%
for 
\begin{eqnarray*}
\widetilde{B} &=& B+U(I-A)^{-1}V,\\
\widetilde{\theta} &=& \theta'+\theta(I-A)^{-1}V,\\
\widetilde{\phi}^{T} &=& \phi'^{T}+U(I-A)^{-1}\phi^{T}.
\end{eqnarray*}
\end{theorem}
This result is related to Theorem~13 of \cite{MT1} for $\beta=1$.
\bigskip

\noindent \textbf{Proof.}
We sketch the proof since it runs along very similar lines to the preceding ones. 
Define a partial permutation graph $\gamma_{\psi}$
for each partial permutation $\psi$ of $\bnp\cup\ms$. 
In this case, the connected subgraphs consist of cycle graphs $\Gamma_{\sigma}$ and  open
necklaces $\Gamma_{\nu}$ containing only $\bn$-vertices, 
and cycles and open necklaces containing at least one $\mathbf{n'}$-vertex.
Define a graph weight $w(\gamma_{\psi})$ as a product of edge weights, vertex weights and cycle factors as before. 
This results in
\begin{equation*}
\sum_{\br}\frac{1}{\br!}
\pperm_{\beta \theta\phi}
\left[\begin{array}{cc}
B & U(\ms) \\ 
V(\ms) & A(\ms,\ms)
\end{array}\right]
=
\frac{e^{\theta(I-A)^{-1}\phi^{T}}}
{\det (I-A)^{\beta }}\sum\limits_{\gamma'\in\,\Gamma'}w(\gamma'),\ 
\end{equation*}
where the sum is over all graphs $\Gamma^{\prime}$ containing at least one $\mathbf{n'}$-vertex. The remaining terms arise as before.

Each $\gamma'\in\Gamma'$ canonically determines a partial permutation $\psi'\in\Psi(\mathbf{n'})$ described by the corresponding ordered sequences of  $\mathbf{n'}$-vertices (for any intermediate $\bn$-vertices). 
As before, the total edge weight coming from all subgraphs joining two $\mathbf{n'}$-vertices 
$i^{\prime }$ and $j^{\prime}$ with intermediate $\bn$-vertices  is 
$\widetilde{B}_{i^{\prime}j^{\prime}}$. 
The total weight arising from the subgraphs of all necklaces joining $\bn$-vertices to
an $\mathbf{n'}$-vertex $i^{\prime}$ with intermediate $\bn$-vertices as depicted in Fig.~5  is  
\begin{eqnarray*}
&&\theta'_{i^{\prime }}
+(\theta V)_{i^{\prime }}
+(\theta A V)_{i^{\prime }}
+\ldots  \\
&=&(\theta'+\theta(I-A)^{-1}V)_{i^{\prime }}=\widetilde{\theta}_{i^{\prime }}.
\end{eqnarray*}%
\begin{center}
\begin{picture}(340,40)


\put(50,20){\makebox(1,0){$i'$}}
\put(50,10){\circle*{4}}
\put(60,10){\makebox(1,0){,}}


\put(80,10){\line(1,0){40}}
\put(80,10){\vector(1,0){23}}
\put(80,20){\makebox(1,0){$j$}}
\put(80,10){\circle*{4}}

\put(120,20){\makebox(1,0){$i'$}}
\put(120,10){\circle*{4}}

\put(130,10){\makebox(1,0){,}}


\put(150,10){\line(1,0){40}}
\put(150,10){\vector(1,0){23}}
\put(150,20){\makebox(1,0){$j$}}
\put(150,10){\circle*{4}}

\put(190,10){\line(1,0){40}}
\put(190,10){\vector(1,0){23}}
\put(190,20){\makebox(1,0){$k$}}
\put(190,10){\circle*{4}}

\put(230,20){\makebox(1,0){$i'$}}
\put(230,10){\circle*{4}}

\put(250,10){\makebox(1,0){$,\ldots$}}

\end{picture}

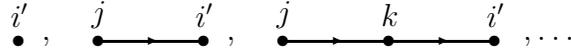
{\small Fig.~5 }
\end{center}
Likewise, the total weight arising from all subgraphs joining an $\mathbf{n'}$-vertex $j^{\prime}$  to
$\bn$-vertices with intermediate $\bn$-vertices is $\widetilde{\phi}_{j^{\prime }}$. 
Combining these results we find that 
\begin{eqnarray*}
\sum\limits_{\gamma'\in\,\Gamma'}w(\gamma') &=&
\pperm_{\beta\,\widetilde{\theta}\,\widetilde{\phi}}\widetilde{B}.\quad \square
\end{eqnarray*}

\bigskip 
\section{The Derangement MMT} 
\label{Sect_DMMT}
Let $\Delta_n\subset\Sigma_{n}$ denote the derangements of the set $\{1,\ldots ,n\}$ 
i.e. each $\pi \in \Delta_{n}$ contains no cycles of length 1.
We introduce the $\beta$-extended Deranged Permanent of an $n\times n$ matrix $A$ by
\begin{equation}
\dperm_{\beta }A=\sum\limits_{\pi \in \Delta_{n}}
\beta^{C(\pi)}\prod\limits_{i}A_{i\pi (i)}.  \label{beta_dperm}
\end{equation}
Using the same multiset notation as before we find
\begin{theorem}
\label{Theorem_DerMMT} 
\begin{equation}
\sum_{\br}
\frac{1}{\br!} \dperm_{\beta} A(\ms,\ms)
=\frac{e^{-\beta \Tr A}}{\det (I-A)^{\beta }}.  \label{dMMT}
\end{equation}
\end{theorem}

\noindent \textbf{Proof.} 
Following the proof of Theorem~\ref{Theorem_MMT} we find
\begin{eqnarray*}
\sum_{\br}
\frac{1}{\br!} \dperm_{\beta} A(\ms,\ms)
&=&
\exp \left( \sum_{\gamma_{\sigma }\in \Gamma_{\sigma },|\sigma|\ge 2}\frac{w(\gamma_{\sigma })}{|\Aut(\gamma_{\sigma })|}\right),
\end{eqnarray*}
where cycles of length one are excluded. Using  
\begin{eqnarray*}
\sum_{\gamma_{\sigma }\in \Gamma_{\sigma },|\sigma|\ge 2}
\frac{w (\gamma_{\sigma })}{|\Aut(\gamma_{\sigma })|} &=&
\beta \sum_{s\geq 1}\frac{1}{s}\Tr(A^{s})-\beta \Tr A \\
&=&-\beta \Tr\log (I-A)-\beta \Tr A,
\end{eqnarray*}
the result follows. $\square$
\bigskip

\textbf{Example.} Consider $n=1$ with $A=z$. Then for multisets $\{1^{r}\}$  we find 
\begin{equation*}
\dperm_{\beta}A(1^{r},1^{r})=d_{r}(\beta)z^{r},
\end{equation*}
where $d_{r}(\beta)=\sum_{s}d_{rs}\beta^{s}$ is the generating polynomial for $d_{rs}$ the 
number of derangements of $r$ labels with $s$ cycles. 
From Theorem~\ref{Theorem_DerMMT} the exponential generating function for $d_{r}(\beta)$ is \cite{HT}
\begin{equation*}
\sum_{r\ge 0} \frac{1}{r!}z^{r}d_{r}(\beta)=
\left( \frac{e^{-z}}{1-z}\right)^\beta.
\end{equation*}
\bigskip

Finally, we can further generalize Theorem~\ref{Theorem_DerMMT} to deranged permanents of submatrices as in Theorem~\ref{Theorem_MMTsub}. Using the notation of \eqref{Block} and \eqref{BAk} we find using similar techniques that
\begin{theorem}
\label{Theorem_dMMTsub} 
\begin{eqnarray}
\sum_{\br}\frac{1}{\br!}\dperm_{\beta }
\left[\begin{array}{cc}
B & U(\ms) \\ 
V(\ms) & A(\ms,\ms)
\end{array}\right] 
&=&\frac{e^{-\beta \Tr A}\cdot \perm_{\beta }\widehat{B}}{\det
(I-A)^{\beta }},
\end{eqnarray}%
for $n'\times n'$ matrix  
\begin{equation*}
\widehat{B}=B-\mathrm{diag}\,B+U(I-A)^{-1}V,
\end{equation*} 
where $\mathrm{diag}\,B_{i'j'}=B_{i'i'}\delta_{i'j'}$.\quad $\square$
\end{theorem}

\end{document}